\documentclass[leqno,11pt]{amsart}

\usepackage{verbatim}
\usepackage{amsmath,amscd,amsthm,amsxtra,amssymb}
\usepackage{epsfig,graphics,color,colortbl}
\usepackage{amssymb,latexsym}
\usepackage{mathrsfs}
\usepackage{soul}   
\usepackage[poly,all]{xy}
\usepackage{marginnote}
\usepackage[colorlinks=true, pdfstartview=FitV, linkcolor=blue,citecolor=blue,urlcolor=blue]{hyperref}

\usepackage[usenames,dvipsnames,svgnames,table]{xcolor}
\usepackage{mathtools}
\def\multiset#1#2{\ensuremath{\left(\kern-.3em\left(\genfrac{}{}{0pt}{}{#1}{#2}\right)\kern-.3em\right)}}

\usepackage[colorlinks=true, pdfstartview=FitV, linkcolor=blue,citecolor=blue, urlcolor=blue]{hyperref}

\newtheorem{thm}{Theorem}[section]
\newtheorem{prop}[thm]{Proposition}
\newtheorem{cor}[thm]{Corollary}             
\newtheorem{lem}[thm]{Lemma}

\theoremstyle{definition}

\newtheorem{Ex}[thm]{Example}

\theoremstyle{remark}
\newtheorem{Rmk}[thm]{Remark}

\newenvironment{red}
{\relax\color{red}}
{\hspace*{.5ex}\relax}

\newcommand{\ber}{\begin{red}}
\newcommand{\er}{\end{red}}

\newcommand{\berE}{\begin{red}{}\marginnote{\fbox{\scshape\lowercase{E}}}{}}
\newcommand{\berY}{\begin{red}{}\marginnote{\fbox{\scshape\lowercase{Y}}}{}}

\numberwithin{equation}{section}

\newcommand{\Z}{\mathbb{Z}}
\newcommand{\Q}{\mathbb{Q}}
\newcommand{\R}{\mathbb{R}}
\newcommand{\C}{\mathbb{C}}

\newcommand{\g}{\mathfrak{g}}
\newcommand{\h}{\mathfrak{h}}

\newcommand{\gl}{\mathfrak{gl}}

\newcommand{\Hom}{\mathrm{Hom}}

\newcommand{\rlQ}{\mathsf{Q}}   
\newcommand{\wlP}{\mathsf{P}}   
\newcommand{\weyl}{\mathsf{W}}  
\newcommand{\cmA}{\mathsf{A}}  
\newcommand{\tf}{\tilde{f}}  
\newcommand{\te}{\tilde{e}}  

\newcommand{\wt}{\mathrm{wt}} 		

\newcommand{\rP}{\Delta}  		

\newcommand{\sg}{\mathfrak{S}}   

\newcommand{\lan}{\langle} 	
\newcommand{\ran}{\rangle}	

\newcommand{\SST}{\mathsf{SST}}					

\newcommand{\La}{\Lambda}
\newcommand{\la}{\lambda}
\newcommand{\al}{\alpha}

\newcommand{\cs}{\mathsf{s}} 			
\newcommand{\ch}{\mathrm{ch}}					
\newcommand{\qdim}{\dim_q}					
\newcommand{\cc}{\mathsf{c}}					
\newcommand{\pr}{\mathsf{pr}}					

\newcommand{\cb}{\mathsf{b}}					
\newcommand{\ca}{\mathsf{a}}
\newcommand{\Orb}{\mathrm{Orb}} 		
\newcommand{\Ula}{T_\la^0}

\begin{document}

\title[$q$-dimensions of highest weight crystals and cyclic sieving phenomenon]
{$q$-dimensions of highest weight crystals and cyclic sieving phenomenon}

\author[Young-Tak Oh]{Young-Tak Oh}
\thanks{The research of Y.-T. Oh was supported by the National Research Foundation of Korea (NRF) Grant funded by the Korean Government (NRF-2020R1F1A1A01071055).}
\address{Department of Mathematics, Sogang University, Seoul 121-742, Republic of Korea}
\email{ytoh@sogang.ac.kr}

\author[Euiyong Park]{Euiyong Park}
\thanks{The research of E. Park was supported by the National Research Foundation of Korea(NRF) Grant funded by the Korea Government(MSIP)(NRF-2020R1F1A1A01065992 and NRF-2020R1A5A1016126).}
\address{Department of Mathematics, University of Seoul, Seoul 02504, Republic of Korea}
\email{epark@uos.ac.kr}

\date{\today}

\begin{abstract}

In this paper, we compute explicitly the $q$-dimensions of highest weight crystals  modulo $q^n-1$ for a quantum group of  arbitrary finite type under certain assumption, and interpret the modulo computations in terms of the cyclic sieving phenomenon.
This interpretation gives an affirmative answer to the conjecture by Alexandersson and Amini.
As an application, under the assumption that $\la$ is a partition of length $<m$
and there exists a fixed point in $\SST_m(\la)$ under the action $\cc$ arising from the crystal structure,
we show that the triple $(\SST_m(\la), \lan \cc \ran, \cs_{\la}(1,q,q^2, \ldots, q^{m-1}))$
exhibits the cycle sieving phenomenon if and only if 
$\la$ is of the form $((am)^{b})$, where either $b=1$ or $m-1$.
Moreover, in this case,  we give an explicit formula to compute the number of all orbits of size $d$
for each divisor $d$ of $n$. 
\end{abstract}

\maketitle


\vskip 2em

\section*{Introduction}

The {\it cyclic sieving phenomenon} was introduced by Reiner-Stanton-White in \cite{RSW04}. 
Let $X$ be  a finite set on which a cyclic group $C$ of order $n$ acts and $f(q) $ a polynomial in $q$ with nonnegative integer coefficients.
We say that $(X, C, f(q))$ exhibits the  cyclic sieving phenomenon if, for all $c\in C$, we have 
$$
\# X^c = f(\zeta_{{\rm o}(c)}),
$$
where $X^c$ is the fixed point set under the action of $c$, 
${\rm o}(c)$ is the order of $c$, and $\zeta_{{\rm o}(c)}=e^{2\pi i/{\rm o}(c)}$.
Many instances of the cyclic sieving phenomenon have been observed for various combinatorial objects including words, multisets, permutations, and tableaux
(see \cite{RSW04,S11} for details).

Let $\pr$ be the {\it promotion} operator due to Sch\"{u}tzenberger \cite{Sch72, Sch77}, and let $\SST_m(\la)$ be the set of semistandard Young tableaux  of shape $\la$ with entries in $\{1,2, \ldots, m\}$.
The cyclic sieving phenomenon about $\SST_m(\la)$ and $\pr$ has drawn a lot of attention from many researchers (see \cite{AOL20, APRU20, BMS14,  BST10, FK14, Rho10, R2001,Wes16} for example).
One of the  most important results in this direction is due to Rhoades \cite{Rho10}, 
who proved that if $\la$ is of rectangular shape, the triple 
$$
\left(\SST_m(\la), ~\langle {\pr} \rangle, ~q^{-\kappa(\la)} s_\la(1,q, \ldots, q^{m-1}) \right) 
$$ 
exhibits  the cyclic sieving phenomenon, where $\kappa(\la)=\sum_{i\ge 1}(i-1)\la_i$,   
and $s_\la(1,q, \ldots, q^{m-1})$ is the principal specialization of the \emph{Schur polynomial} $s_\la(x_1, x_2, \ldots, x_m)$.
However, this result is no longer valid outside rectangular shape in general.
If one wants to keep $\SST_m(\la)$ and the principal specialization of the Schur polynomial in the triple, another appropriate action other than $\pr$ should be considered.
In a previous article of the authors, \cite{OP19}, a new cyclic sieving phenomenon triple 
\begin{align} \label{Eq: new triple}
\left(\SST_m(\la), ~\langle {\cc} \rangle, ~q^{-\kappa(\la)} s_\la(1,q, \ldots, q^{m-1}) \right) 
\end{align}
was provided
under the condition $\gcd(|\la|, m)=1$, where the action $\cc$ arises \emph{naturally} from the $U_q(\mathfrak{sl}_m)$-\emph{crystal} structure of $\SST_m(\la)$. 
Crystal bases theory is one of the most powerful combinatorial tools for studying representations of quantum groups in the viewpoint of graph theory with natural connections to tableaux and functions invariant under the action of the Weyl group like symmetric functions (\cite{BS17, HK02, Kas90, Kas91, Kas93}).
The promotion operator $\pr$ on $\SST_m(\la)$ with entries $\le m$ can be defined by $\pr := \sigma_1 \sigma_2 \cdots \sigma_{m-1}$, where the $\sigma_i$ are the \emph{Bender-Knuth involutions}, certain natural involutions on tableaux that exchange the number of $i$’s and $i+1$’s. 
The operator $\cc$ on $\SST_m(\la)$ can similarly be defined by $\cc := s_1 s_2 \cdots s_{m-1}$, where the $s_i$ are the generators of the symmetric group $\sg_m$ action on the crystal $\SST_m(\la)$. 
Since $s_1 s_2 \cdots s_{m-1}$ is a Coxeter element of $\sg_m$, it has order $m$. Thus the cyclic group $C$ of order $m$ acts on $\SST_m(\la)$ via the operator $\cc$; for rectangular shape partitions, the action of $\pr$ has order $m$, but for other shapes it does not  because the Bender-Knuth involutions do not give an action of the symmetric group.
We remark that, before crystal theory was developed, the same symmetric group action on $\SST_m(\la)$ was studied at a purely combinatorial level by Lascoux and Sch\"{u}tzenberger in \cite{LS78}.

Without the condition $\gcd(|\la|, m)=1$, the new triple $\eqref{Eq: new triple}$  does not exhibit the cyclic sieving phenomenon in general. 
Thus, it is an interesting problem to find a necessary and sufficient condition for the  cyclic sieving phenomenon of the new triple $\eqref{Eq: new triple}$.
To answer this problem, we follow the method of Alexandersson and Amini \cite{AA19}.
To be precise, we ask what conditions guarantee the existence
of an action of a cyclic group $C$ of order $n$ on $\SST_m(\la)$, without being able to describe it explicitly, 
such that the triple $\left(\SST_m(\la), C, q^{-\kappa(\la)} s_\la(1,q, \ldots, q^{m-1}) \right) $ exhibits the cyclic sieving phenomenon.

In this paper, we compute explicitly the \emph{$q$-dimensions} of highest weight crystals  modulo $q^n-1$ for a quantum group of \emph{arbitrary finite type} under certain assumptions, and interpret the modulo computations in terms of the cyclic sieving phenomenon.
Let $\g$ be a finite-dimensional simple Lie algebra over $\C$ and $U_q(\g)$ be its quantum group. We write $\rP^+$ for the set of positive roots of $\g$. For a dominant integral weight $\La$, let $B(\La)$ be the highest weight $U_q(\g)$-crystal with highest weight $\La$.
We denote by $\dim_q B(\La)$ the $q$-dimension of $B(\La)$, which is the polynomial in $q$ obtained from the character $\ch B(\La)$ by specializing at $q^{(\rho,-)}$ (see Section \ref{Sec: Preliminaries} for the definition).
When $U_q(\g)$ is of type $A_{m-1}$, i.e., $\g = \mathfrak{sl}_m$, the crystal $B(\La)$ can be realized as $\SST_m(\la)$ and the $q$-dimension $\dim_q B(\La)$ is equal to the principal specialization $q^{-\kappa(\la)} s_\la(1,q, \ldots, q^{m-1})$ of the Schur polynomial $s_\la(x_1,x_2, \ldots, x_m)$.
Here, $\La$ and $\la$ are related in $\eqref{eq: la La}$. 
Let $n$ be a positive integer. Under the assumption that 
\begin{align} \label{Eq: Con_int}
\text{ $(\beta, \La)$ is divisible by $n$ for any $\beta \in \rP^+$,} 
\end{align}
we provide an explicit expression for $\qdim B(\La)$  modulo $q^n-1$ using the Weyl character formula as follows:  
$$
\qdim B(\La) \equiv \sum_{d|n} \ca_d \,\, \frac{q^n-1}{q^{\frac nd}-1} \pmod{q^n-1},
$$
where $\ca_d \in \Z_{\ge0} $ are nonnegative integers given explicitly in \eqref{Eq: def of a b}
(see Theorem~\ref{Thm: main}).
Note that when $U_q(\g)$ is of type $A_{m-1}$ and $n=m$, 
condition \eqref{Eq: Con_int} implies $|\lambda|$ being divisible by $m$, i.e., $\gcd(|\la|, m) =m$, 
which case is not covered by the previous result of the authors, \cite[Theorem 4.3]{OP19}.
We can also derive a similar result for the $q$-dimension
$\dim_q^{\vee} B(\La)$ of $B(\La)$, which is obtained by specializing at $q^{\lan \rho^{\vee},-\ran}$ (see Remark \ref{second congruence}).
It should be remarked that there are root of unity evaluations of $\dim_q^{\vee} B(\La)$ that have been studied in the literature.
For instance, letting $\varphi_\La(q):=q^{-\lan \rho^{\vee}, \Lambda\ran}\qdim^{\vee}  B(\La)$, it is known by Kac (\cite[Exercise 10.15]{Kac81} or \cite{Kac90}) that if $\omega$ is a root of unity of order equal to the Coxeter number of the Weyl group, then $\varphi_\La(\omega)=0$ or $\pm 1$.

From the viewpoint of the cyclic sieving phenomenon, 
the above computation modulo $q^n-1$ says that there exists an action of a cyclic group $C$ of order $n$ 
on $B(\La)$, 
 without being able to describe it explicitly,
such that 
the triple $(B(\La), C, \qdim B(\La))$ exhibits the cyclic sieving phenomenon
and the number of all orbits of size $d $ is equal to $\ca_d$ for any positive integer $d$ with $d | n$ (see Theorem~\ref{Thm: csp and B}).
In the case where $\g = \gl_m $, the situation is more interesting. Let $\la = (\la_1 \ldots, \la_{\ell})$ be a partition and $\La$ the dominant integral weight given in $\eqref{eq: la La}$. 
In this case, the condition $\eqref{Eq: Con_int}$ is equivalent to 
$$
\text{ $\la_i - \la_j$ is divisible by $n$ for all $1 \le i < j \le m$,}
$$
which means that $\la$ is a \emph{stretched} Young diagram by $n$, i.e., $\la = n \tilde{\la}$  for some Young diagram $\tilde{\la}$.
Hence Theorem~\ref{Thm: csp and B} implies that, for a stretched Young diagram $\la$ by $n$,
there exists an action of a cyclic group $C$ of order $n$ on $\SST_m(\la)$ such that 
the triple 
$$
(\SST_m(\la),  C , q^{-\kappa(\la)} s_\la( 1,q, q^2, \ldots, q^{m-1})  )
$$
exhibits the cyclic sieving phenomenon. Consequently, we give an affirmative answer to the conjecture \cite[Conjecture 3.4]{AA19} by Alexandersson and Amini (see Corollary \ref{criterion in type A}). 
In this viewpoint, Theorem~\ref{Thm: csp and B} can be understood as an affirmative answer to a crystal-theoretical generalization of 
this conjecture.

We next focus on the case where $\mathfrak g$ is of type $A_{m-1}$ and 
the action $\cc$ arising from the crystal structure.
In the previous article of the authors, \cite{OP19},
the case where $\gcd(m, |\la|)=1$ was studied extensively, 
where every orbit is free.  
We now consider the case at least one fixed point exists. 
Under the assumption that $\la$ is a partition of length $<m$
and there exists a fixed point in $\SST_m(\la)$ under the action of $\cc$,
we show that the triple $(\SST_m(\la), \lan \cc \ran, \cs_{\la}(1,q,q^2, \ldots, q^{m-1}))$
exhibits the cycle sieving phenomenon if and only if 
$\la$ is of the form $((am)^{b})$, where either $b=1$ or $m-1$
(Theorem~\ref{non relatively prime case}). Moreover, in this case,  we give an explicit formula to compute the number of all orbits of size $d $ (Proposition \ref{Prop: formula for orbits}).
When $m$ is a prime $p$, 
combining this result with \cite[Theorem 4.3]{OP19} enables us to characterize  
when the triple $(\SST_p(\la), \lan \cc \ran, \cs_{\la}(1,q,q^2, \ldots, q^{p-1}))$
exhibits the cycle sieving phenomenon.
Because there exists a fixed point in $\SST_m(\la)$ under the action of $\cc$ if and only if $|\lambda|$ is divisible by $m$, the problem of when $(\SST_p(\la), \lan \cc \ran, \cs_{\la}(1,q,q^2, \ldots, q^{p-1}))$ exhibits the cyclic sieving phenomenon will be completely settled if one can successfully attack the cases $1< \gcd(m, |\la|)<m$.

This paper is organized as follows.
In Section~\ref{Sec: Preliminaries}, we introduce the prerequisites 
on highest weight crystals and their $q$-dimensions.
In Section~\ref{Congruencs of the $q$-dimension},
we derive a congruence relation of $\qdim B(\La)$
which plays a crucial role throughout this paper. 
In Section~\ref{Cyclic sieving phenomena and $q$-dimensions}, 
we reinterpret the congruence obtained in Section~\ref{Congruencs of the $q$-dimension} in the viewpoint of 
the cyclic sieving phenomenon
and apply it to the case where $\mathfrak g =\mathfrak{gl}_m$.
The final section is devoted to characterizing   
when the triple $(\SST_p(\la), \lan \cc \ran, \cs_{\la}(1,q,q^2, \ldots, q^{p-1}))$
exhibits the cycle sieving phenomenon under the assumption that there exists at least one fixed point.

\vskip 2em

\noindent
{\bf Acknowledgments.}
The authors would like to thank the anonymous reviewers for their valuable comments and suggestions.

\vskip 2em
\section{Highest weight crystals and their $q$-dimensions} \label{Sec: Preliminaries}

Let $I$ be a finite index set and let $\cmA = (a_{ij})_{i,j\in I}$ be a \emph{Cartan matrix} of \emph{finite} type.
We choose a diagonal matrix $D={\rm diag}(\mathsf d_i \in \Z_{\ge 0} \mid i \in I)$ such that $\min\{ \mathsf{d}_i \mid i\in I \}=1$ and 
$D\cmA$ is symmetric.
We then consider a quintuple $ (\cmA,\wlP,\Pi,\wlP^\vee,\Pi^\vee) $, called {\it a Cartan datum associated with $A$}, such that  
\begin{enumerate}
\item $\wlP$ is a free abelian group of rank $|I|$, called the {\em weight lattice},
\item $\Pi = \{ \alpha_i \mid i\in I \} \subset \wlP$,
called the set of {\em simple roots},
\item $\wlP^{\vee}=
\Hom_{\Z}( \wlP, \Z )$, called the \emph{coweight lattice}, 
\item $\Pi^{\vee} =\{ h_i \in \wlP^\vee \mid i\in I\}$, called the set of {\em simple coroots},
\end{enumerate}
which satisfy the following requirements: 
\begin{itemize}
\item $\lan h_i, \alpha_j \ran = a_{ij}$ for $i,j \in I$,
\item $\Pi$ is linearly independent over $\Q$, and 
\item for each $i\in I$, there exists $\varpi_i \in \wlP$, called the \emph{fundamental weight}, such that $\lan h_j,\varpi_i \ran =\delta_{j,i}$ for all $j \in I$.
\end{itemize}
We denote by $\wlP^+: =\{ \lambda \in \wlP \mid \lan h_i, \lambda\ran \ge 0 \ \text{for all }  \ i \in I \}$ the set of \emph{dominant integral weights}.
There exists a nondegenerate symmetric bilinear form $( \cdot \, , \cdot )$ on $ \wlP$ satisfying
\begin{equation*}
(\alpha_i,\alpha_j)=\mathsf d_i a_{ij} \quad (i,j \in I),\quad \text{and } \quad  \lan h_i,  \lambda\ran = \dfrac{2 (\alpha_i,\lambda)}{(\alpha_i,\alpha_i)} \quad (\lambda \in \wlP, \ i \in I).
\end{equation*}
Let $ \rlQ := \bigoplus_{i \in I} \Z \alpha_i$ be the \emph{root lattice}, and
let $\rP \subset \rlQ$ be the set of \emph{roots} associated with $\cmA$. We write $\rP^+$ for the set of positive roots.

Fix an indeterminate $q$.
Let $U_q(\g)$ be the \emph{quantum group} associated with $(\cmA, \wlP,\wlP^\vee \Pi, \Pi^{\vee})$, which
is the associative algebra over $\mathbb C(q)$ with $1$ generated by $f_i$, $e_i$ $(i\in I)$ and $q^h$ $(h\in \wlP)$ with certain defining relations (see \cite[$\S$3]{HK02} for details).
For a dominant integral weight $\La \in \wlP^+$, we denote by $V_q(\La)$ the irreducible highest weight $U_q(\g)$-module with highest weight $\La$, and 
denote by $B(\La)$ its crystal. 
We denote by $ \te_i$ and $\tf_i$ ($i\in I$) the crystal operators on $B(\La)$. 
We refer the reader to \cite{ BS17, HK02,Kas90, Kas91, Kas93} for crystals.

We set $\weyl$ to be the \emph{Weyl group} associated with $\cmA$, which is a subgroup of $\mathrm{Aut}(\wlP)$ generated by
$s_i(\la) := \la - \langle h_i, \la \rangle \alpha_i$ for $i\in I$ and $ \la\in \wlP$.
Note that $\rP$ is invariant under the actions of $\weyl$. The Weyl group $\weyl$ also acts on the crystal $B(\La)$ as follows:
for $i\in I$ and $b\in B(\La)$, define 
\begin{align} \label{Eq: def of si}
s_i ( b) :=
\left\{
\begin{array}{ll}
 \tf_i^{ \langle h_i, \wt(b) \rangle } b  & \text{ if }  \langle h_i, \wt(b) \rangle \ge 0,\\
\te_i^{ - \langle h_i, \wt(b) \rangle }b & \text{ if } \langle h_i, \wt(b) \rangle < 0.
\end{array}
\right.
\end{align}

We set $B(\La)_\xi := \{ b\in B(\La) \mid \wt(b) = \xi \}$ so that $B(\La) = \sqcup_{\xi \in \wlP} B(\La)_\xi$, and $\wt(B(\La)) = \{ \mu \in \wlP \mid B(\La)_\mu \ne \emptyset                                                                                                                                               \}$.
The \emph{character} $\ch B(\La)$ of $B(\La)$ is defined by
$$
\ch B(\La) := \sum_{\xi \in \wt(B(\La))}  | B(\La)_\xi | e^\xi,
$$
where $| B(\La)_\xi |$ is the number of elements of $B(\La)_\xi$, and
$e^\xi$ are formal basis elements of the group algebra $\Q[\wlP]$ with the multiplication given by $e^\xi e^{\xi'} = e^{\xi+\xi'}$.
The \emph{Weyl character formula} says that 
$$
\ch B(\La) = \frac{ \sum_{w\in \weyl} (-1)^{\ell(w)} e^{w( \La+\rho) - \rho   }  }{ \prod_{\beta\in \rP^+} (1-e^{- \beta} )},
$$
where $\rho = \frac{1}{2}\sum_{\beta \in \rP^+} \beta$ (for instance, see \cite[Theorem 10.4]{Kac90}). 
Note that $\rho =  \sum_{i\in I} \varpi_i$.

In this paper, we consider the following polynomials in $q$ arising from the crystal $B(\La)$: 
\begin{equation*}
\begin{aligned}
&\qdim B(\La):=\sum_{\xi \in \wt(B(\La))}  | B(\La)_\xi | q^{ ( \rho, \La - \xi )},\\	
&\qdim^{\vee} B(\La):=\sum_{\xi \in \wt(B(\La))}  | B(\La)_\xi | q^{\lan \rho^{\vee}, \La - \xi \ran},
\end{aligned}
\end{equation*}
where $  \vee : \h^* := \C \otimes \wlP \longrightarrow  \h:= \C \otimes \wlP^\vee $ is the isomorphism given in \cite[$\S$2.1]{Kac90}.  
These polynomials can be obtained via specializations of $\ch B(\La)$.
Define homomorphisms 
\begin{equation*}
\begin{aligned}
&F_\rho: \mathbb C[\rlQ]\to \mathbb  C[q^{\pm1}], \quad e^\lambda \mapsto q^{-(\rho, \lambda)},\\
&F_{\rho^{\vee}}: \mathbb C[\rlQ ]\to \mathbb C[q^{\pm1}], \quad e^\lambda \mapsto q^{- \lan \rho^{\vee}, \lambda \ran }.
\end{aligned}
\end{equation*}	
One easily sees that $\qdim B(\La)=F_\rho(e^{-\Lambda} \ch B (\Lambda))$ and $\qdim^\vee B(\La)=F_{\rho^{\vee}}(e^{-\Lambda} \ch B(\Lambda))$.
Using the Weyl character formula and properties of $\rho$, one can show that 
\begin{align}
& \label{Eq: q-dim}\qdim B(\La)=  \prod_{\beta\in \rP^+} \frac {1-q^{(\beta, \Lambda+\rho)}}{1-q^{(\beta, \rho)}}, \\
& \label{Eq: second q-dim}\qdim^{\vee}  B(\La)=\prod_{\beta\in \rP^+} \frac {1-q^{\lan \beta^{\vee}, \Lambda+\rho \ran}}{1-q^{\lan \beta^{\vee}, \rho \ran}}.
\end{align}
(see \cite[$\S$10.10]{Kac90} and \cite{Stem94,Stem96}).
In the literature such as \cite[$\S$10.10]{Kac90}, the right hand side of \eqref{Eq: second q-dim} is 
called {\it the $q$-dimension of $V(\La)$}.
Note that $ \qdim B(\La) = \qdim^\vee B(\La)$ if the Cartan matrix $\cmA$ is symmetric.

\vskip 2em

\section{Congruencs of the $q$-dimension of $B(\La)$} \label{Congruencs of the $q$-dimension}

In this section, we provide some noteworthy congruences of the $q$-dimensions of $B(\La)$ which are significant not only in itself but also 
a cyclic sieving phenomenon on $B(\La)$.
We fix a Cartan matrix $\cmA$ of finite type and its Cartan datum $ (\cmA,\wlP,\Pi,\wlP^\vee,\Pi^\vee) $.
Let $\La \in \wlP^+$ and $n\in \Z_{>0}$. 
From now on, we impose the following condition:
\begin{align} \label{Assumption}
\text{ {\bf{Condition}} :  for any $\beta \in \rP^+$, $(\beta, \La)$ is divisible by $n$.}
\end{align}
For any $d | n$, we set 
$$
\rP^+_d := \{ \beta \in \rP^+ \mid \text{ $(\beta, \rho)$ is divisible by $d$} \}. 
$$
Denote by $\Phi_d(q)$ the $d$th \emph{cyclotomic polynomial}. Note that $\Phi_d(q)$ is irreducible over $\Z$ and 
\[
q^n-1= \prod_{d|n}\Phi_d(q).
\]

\begin{lem} \label{Lem: qdim modulo} 
Assume that the condition \eqref{Assumption} holds.
For $d \mid n$, we have 
$$
\qdim B(\La) \equiv \prod_{\beta \in \rP^+_d} \left(\frac{ (\beta, \La )}{ (\beta, \rho)} +1 \right) \pmod{\Phi_d(q)}.
$$
Here, the right hand side is set to be $1$ in the case where $\rP^+_d=\emptyset$. 
\end{lem}
\begin{proof}
Note that 
\begin{enumerate}
\item[(i)] $q^{(\beta, \La)} \equiv 1 \pmod{\Phi_d(q)}$ for $\beta \in \rP^+$, and 
\item[(ii)] $1-q^{(\beta, \rho)} \not\equiv 0 \pmod{\Phi_d(q)}$ for $\beta \in \rP^+ \setminus \rP^+_d $.
\end{enumerate}
It follows from $\eqref{Eq: q-dim}$ that 
\begin{align*}
\qdim B(\La) & = \left( \prod_{\beta \in \rP^+_d}  \frac{  1- q^{(\beta, \La+\rho )} }{ 1- q^{(\beta, \rho )} } \right) \left( \prod_{\beta \in \rP^+ \setminus \rP^+_d }  \frac{  1- q^{(\beta, \La )} q^{(\beta, \rho )} }{ 1- q^{(\beta, \rho )} } \right)\\
&\equiv  \prod_{\beta \in \rP^+_d}  \frac{  1- q^{(\beta, \La+\rho )} }{ 1- q^{(\beta, \rho )} } 
\equiv  \prod_{\beta \in \rP^+_d}  \frac{  (1-q^d)  (1+q^d + q^{2d} + \cdots + q^{ d \left( \frac{(\beta, \La+\rho)}{d}-1 \right)} )  }{ (1-q^d)  (1+q^d + q^{2d} + \cdots + q^{     d \left(  \frac{(\beta, \rho)}{d}-1 \right)} )  } \\
&\equiv  \prod_{\beta \in \rP^+_d} \frac{ (\beta, \La+\rho) }{ (\beta, \rho) } \pmod{\Phi_d(q)},
\end{align*}
which completes the proof.
\end{proof}

For $d \mid n$, we set 
\begin{equation}
\begin{aligned} \label{Eq: def of a b}
&\cb_{d} :=  \prod_{\beta \in \rP^+_{n/d}} \left(\frac{ (\beta, \La )}{ (\beta, \rho)} +1 \right),\\
&\ca_{d} := \frac{1}{d} \sum_{e \mid d} \mu \left(\frac{d}{e} \right) \cb_e, 
\end{aligned}
\end{equation}
where $\mu$ is the classical M\"{o}bius function. 
Note that 
\begin{align} \label{Eq: Mobius}
\cb_{d} = \sum_{ e\mid d} e \ca_e.
\end{align}
\begin{Ex}\label{Prop: a b} \
From the definition, it is easy to see the following.
\begin{enumerate}
\item Let $\beta_0 \in \rP^+$ be the highest root.
If $ (\beta_0, \rho) < n$, then $\rP^+_n=\emptyset$ and therefore $\ca_1=\cb_1=1$.
\item Let $n=p^l$ for  a prime $p$.
For $0 < s  \le l$, we have
$$
\ca_{p^s} = \frac{1}{p^s} (\cb_{p^s} - \cb_{p^{s-1}} ). 
$$
\end{enumerate}
\end{Ex}

\begin{lem}\label{nonnegarivity lemma}
Let $S$ be a finite set and let $\xi: S \rightarrow \Z_{> 0}$ be a function.
For $s \in S$, let $x_s$ be an indeterminate. For $d | n$, let
\begin{align*}
B_d = \prod_{s \in S, \ \frac{n}{d} | \xi(s)} \left( 1+x_s \right) \quad \text{ and } \quad
A_d = \sum_{k | d} \mu(k) B_{\frac{d}{k}}.  	
\end{align*}
Then we have
$$
A_d = \sum_{T \subset S, \ \gcd( n, \xi(T) )= \frac{n}{d}} x_T,
$$
where $ \displaystyle x_T = \prod_{t\in T} x_t$ and $\xi(T) = \gcd \{ \xi(t) \mid t \in T \}$ for any $T \subset S$.
Here we set $ \xi(\emptyset)=0$.
\end{lem}
\begin{proof}
It follows from
$$
B_d = \prod_{s \in S, \ \frac{n}{d} | \xi(s)} \left( 1+x_s \right) = \sum_{T \subset S, \ \frac{n}{d} | \xi(T)} x_T
$$	
that
\begin{align*}
A_d &= \sum_{k | d} \mu(k) \left( \sum_{T \subset S, \ \frac{nk}{d} | \xi(T)} x_T \right) = \sum_{T \subset S} x_T \left( \sum_{k | d, \ \frac{nk}{d} | \xi(T)} \mu(k) \right) \\
&= \sum_{T \subset S} x_T \left( \sum_{ \frac{nk}{d} | \gcd \left( n, \xi(T) \right)} \mu(k) \right).
\end{align*}	
Thus the assertion follows from
$$\displaystyle
\sum_{ \frac{nk}{d} | \gcd \left( n, \xi(T) \right)} \mu(k)
= \begin{cases}
	1 & \text{ if }  \gcd \left( n, \xi(T) \right) = \frac{n}{d}, \\
	0 & \text{otherwise.}
\end{cases}
 $$
\end{proof}

\begin{Rmk}
In the first arXiv version of the paper, our original lemma used to prove Theorem \ref{Thm: main}
is different from Lemma \ref{nonnegarivity lemma}.
We are very grateful to the anonymous referee for informing us of the present lemma. 
We also would like to thank Martin Rubey for letting us know the same idea independently.
\end{Rmk}

We now shall prove the following congruence for $\qdim B(\La)$.

\begin{thm} \label{Thm: main} \
Assume that $(\beta, \La)$ is divisible by $n$ for any $\beta \in \rP^+$.
\begin{enumerate}
\item For any $d | n$,   $\ca_d \in \Z_{\ge0}$. 
\item 
$$
\qdim B(\La) \equiv \sum_{d|n} \ca_d \,\, \frac{q^n-1}{q^{\frac nd}-1} \pmod{q^n-1}.
$$
\end{enumerate}
\end{thm}
\begin{proof}
We set
$$
B(q):= \sum_{e|n} \ca_e \,\, \frac{q^n-1}{q^{\frac{n}{e}}-1}.
$$

Let $d \in \Z_{>0} $ with $d | n$. 
Since $ \frac{q^n-1}{q^{\frac{n}{e}}-1} \equiv 0 \pmod{\Phi_d(q)}$ for $ d \nmid \frac{n}{e} $ and 
$ \frac{q^n-1}{q^{\frac{n}{e}}-1} \equiv e \pmod{\Phi_d(q)}$ for $ d \mid \frac{n}{e} $, 
by $\eqref{Eq: Mobius}$, we have 
$$B(q)\equiv  \sum_{e|\frac nd}e \ca_e = \cb_{\frac nd} \pmod {\Phi_d(q)}. $$
Combining Lemma \ref{Lem: qdim modulo} with the Chinese Remainder Theorem, we conclude that
\[ \qdim B(\La) \equiv B(q) \pmod {q^n -1}.\]

To complete the proof, it remains to see that $\ca_d \in \mathbb Z_{\ge 0}$ for all $d|n$.
Since $\frac{ (\beta, \La )}{ (\beta, \rho)} \ge 0$ for all $\beta \in \rP^+$,
the non-negativity of $\ca_d$ follows by applying Lemma \ref{nonnegarivity lemma} to the setting $S=\rP^+$ and $\xi(\beta) = (\beta, \rho)$.
We now see that $\ca_d \in \mathbb Z$.
Note that $B(q)$ is the remainder of $\qdim B(\La)$ when divided by $q^n-1$.
Since 
$$
\qdim B(\La) \in \mathbb Z[q]
$$
 and $q^n-1$ is monic, it follows that
$B(q) \in \mathbb Z[q].$
Let $d_1> d_2 \cdots > d_s$ be all divisors of $n$ such that $\ca_{d} \ne 0$.
Then the leading coefficient of $B(q)$ is $\ca_{d_1}$, so it is an integer.
Next, consider $$B(q)-\ca_{d_1}(1+q^{\frac{n}{d_1}}+\cdots + (q^{\frac{n}{d_1}})^{d_1 -1}) \in \mathbb Z[q].$$
Its leading coefficient is given by $\ca_{d_2}-\ca_{d_1}$ if $d_2| d_1$ and $\ca_{d_2}$ otherwise, thus $\ca_{d_2} \in \mathbb Z$.
In this way, we can see inductively that $\ca_d \in \mathbb Z$ for all $d|n$.
\end{proof}

\begin{Rmk}\label{second congruence}
We can also derive an analogue of Theorem \ref{Thm: main} for $\qdim^{\vee} B(\La)$ in the same manner.
For any $d | n$, we set 
\begin{equation}
\begin{aligned}\label{Eq: def of a b II}
&({\rP}^{+}_d)^{\vee} := \{ \beta \in \rP^+ \mid \text{ $\lan \beta^{\vee}, \rho \ran$ is divisible by $d$}\}, \\
& \cb^{\vee}_{d} := \prod_{\beta \in ({\rP}^{+}_{n/d})^{\vee}} \left(\frac{ \lan \beta^{\vee}, \La \ran}{ \lan \beta^{\vee}, \rho \ran} +1 \right),\\ 
&\ca^{\vee}_{d} := \frac{1}{d} \sum_{e \mid d} \mu \left(\frac{d}{e} \right) \cb^{\vee}_e.
\end{aligned}
\end{equation}
Assume that $\lan \beta^{\vee}, \La \ran$ is divisible by $n$ for all $\beta \in \rP^+$. 
In the same manner as above, we can derive that  
\begin{enumerate}
\item For any $d | n$,   $\ca^{\vee}_d \in \Z_{\ge0}$. 
\item 
$$
\qdim^{\vee} B(\La) \equiv \sum_{d|n} \ca^{\vee}_d \,\, \frac{q^n-1}{q^{\frac nd}-1} \pmod{q^n-1}.
$$
\end{enumerate}
\end{Rmk}

\vskip 2em

\section{Cyclic sieving phenomena and $q$-dimensions} \label{Cyclic sieving phenomena and $q$-dimensions} 

As before, let $(\cmA,\wlP,\Pi,\wlP^\vee,\Pi^\vee)$ be a Cartan datum of finite type.
We here interpret Theorem \ref{Thm: main} from the viewpoint of the cyclic sieving phenomenon.
To do this, we need the following lemma.

\begin{lem} {\rm (Alexandersson and Amini~\cite[Theorem 2.7]{AA19})}\label{necessary and sufficient condition for csp}
Let $f(q)\in \mathbb N[q]$ and suppose $f(\omega_n^j)\in \mathbb N$ for each $j=1, \ldots, n$, where $\omega_n$
denotes a primitive $n$th root of unity. 
Let $X$ be any set of size $f(1)$. 
Then there exists an action of a cyclic group $C$ of order $n$ on $X$ such that $(X, C, f(q))$
exhibits the cyclic sieving phenomenon if and only if for each $k|n$,
\[\sum_{j|k}\mu(k/j)f(\omega_n^j)\ge 0.\]
Here, $\mu$ is the M\"{o}bius function.
\end{lem}

Let $X$ be a finite set with an action of a finite cyclic group $C$. 
We denote by $\Orb_C^d(X)$ the set of all orbits of size $d$.
With this notation, we state the following theorem.

\begin{thm} \label{Thm: csp and B}
Let $\La\in \wlP^+$ and $n\in \Z_{>0}$. 
Assume that $(\beta, \La)$ is divisible by $n$ for any $\beta \in \rP^+$.
\begin{enumerate}
\item There exists an action of a cyclic group $C$ of order $n$ on $B(\La)$ such that
the triple $(B(\La), C, \qdim B(\La))$ exhibits the cyclic sieving phenomenon.
\item Let $C$ be a cyclic group of order $n$ acting on $B(\La)$.
\begin{enumerate}
\item  The triple $(B(\La), C, \qdim B(\La))$ exhibits the cyclic sieving phenomenon if and only if
the number of orbits in $\Orb_C^d(B(\La))$ is equal to $\ca_{d}$ for all $d | n$,
where $\ca_{d}$ is given in \eqref{Eq: def of a b}.
\item In particular, if $n$ is prime and $(\beta_0, \rho) <n$, then
the triple $(B(\La), C, \qdim B(\La))$ exhibits the cyclic sieving phenomenon if and only if
$C$ acting on $B(\La)$ has exactly one fixed point.
Here, $\beta_0$ is the highest root and $\rho = \frac{1}{2} \sum_{\beta \in \Delta^+} \beta $.
\end{enumerate}
\end{enumerate}
\end{thm}

\begin{proof}
(1) Let
$$f(q)=\sum_{d|n} \ca_d \,\, \frac{q^n-1}{q^{\frac nd}-1},$$
where $\ca_d$'s are given in \eqref{Eq: def of a b}.
For a divisor $k$ of $n$, one can easily see that 
\[\frac{(\omega_n^k)^n-1}{(\omega_n^k)^{\frac nd}-1}
=\begin{cases}
d & \text{ if } d|k,\\
0 & \text{ if } d \nmid k.
\end{cases}
\]
Combining this with the fact $f(\omega_n^j)=f(\omega_n^{\gcd(n,j)})$, one can check that 
$f(\omega_n^j)\in \mathbb N$ for each $j=1, \ldots, n$.  
For a divisor $k$ of $n$, note that  
$$f(\omega_n^k)=\sum_{d|k}da_d.$$
Applying the M\"{o}bius inversion formula to this equality implies that   
\begin{equation}\label{number of orbits of size k}
\sum_{j|k}\mu(k/j)f(\omega_n^j)=ka_k \ge 0.
\end{equation}
Therefore, the assertion is obtained by combining Theorem \ref{Thm: main} with Lemma \ref{necessary and sufficient condition for csp}.

(2) The assertion (a) follows from Theorem \ref{Thm: main} and \cite[Proposition 2.1]{RSW04}
and the assertion (b) follows from Example \ref{Prop: a b}.
\end{proof}

\begin{Rmk}
Assume that $\lan \beta^{\vee}, \La \ran$ is divisible by $n$ for all $\beta \in \rP^+$.
Due to Remark~\ref{second congruence}, 
we can also derive an analogue of Theorem \ref{Thm: csp and B} for the triple 
$(B(\La), C, \qdim^{\vee} B(\La))$.
\end{Rmk}

Let us consider the case where $U_q(\g)=U_q(\mathfrak{gl}_m)$.
Let $I=\{1,2, \ldots,m-1 \}$ and let $\{ \epsilon_1, \ldots, \epsilon_m \}$ be the standard orthonormal basis of the Euclidean space $\R^m$. Then one can realize the weight lattice of $U_q(\mathfrak{gl}_m)$ inside $\R^m$ (see \cite[$\S$2.1]{BS17} and \cite[$\S$7.1]{HK02} for example).
In this realization, we have 
$$
\rP_+ = \{ \epsilon_i - \epsilon_j \mid 1 \le i < j \le m \}.
$$
Let $\la = (\la_1, \la_2, \ldots, \la_{\ell})$ be a partition with $\ell \le m$ and set 
\begin{align} \label{eq: la La}
\La :=  \sum_{k=1}^{\ell} \la_k \epsilon_k
\end{align}
Note that 
$$
 \qdim B(\La)= q^{-\kappa(\la)}s_\la( 1,q, q^2, \ldots, q^{m-1}), 
$$
where 
$\la= (\la_1, \la_2, \ldots)$ is the Young diagram of length $\le m$ corresponding to $\La$, $\kappa(\la) = \sum_{ k \ge 1} (k-1)\la_k$, and 
$s_\la(1,q, \ldots, q^{m-1})$ is the principal specialization of the \emph{Schur polynomial} $s_\la(x_1, x_2, \ldots, x_{m})$.

It is straightforward to show that 
$(\beta, \La)$ is divisible by $n$ for all $\beta \in \rP^+$ if and only if 
\begin{align} \label{Eq: stretched}
\text{ $\la_i - \la_j$ is divisible by $n$ for all $1 \le i < j \le m$. }
\end{align}
Moreover, we know that 
\begin{enumerate}
\item[(a)]  the set $\SST_m(\la)$ of all semistandard tableaux of shape $\la$ with entries $\{1,2, \ldots, m \}$ has a $U_q(\mathfrak{gl}_m)$-crystal structure which is isomorphic to $B(\La)$, and 
\item[(b)] $\qdim B(\La)  = q^{-\kappa(\la)} s_\la( 1,q, q^2, \ldots, q^{m-1}) $.
\end{enumerate} 
Applying these facts to Theorem \ref{Thm: csp and B} yields the following corollary,  
which gives an affirmative answer for the conjecture in \cite[Conjecture 3.4]{AA19}. 

\begin{cor} \label{criterion in type A}
Let $m, n \in \Z_{>0}$ and let $\la=(\la_1, \la_2, \ldots, \la_{\ell} )$ be a partition such that  
$\la_i - \la_j$ is divisible by $n$ for all $1 \le i < j \le m$. 
Then there exists an action of a cyclic group $C$ of order $n$ on $\SST_m(\la)$ such that 
the triple $(\SST_m(\la),  C , q^{-\kappa(\la)} s_\la( 1,q, q^2, \ldots, q^{m-1})  )$ exhibits the cyclic sieving phenomenon.
\end{cor}

\begin{Rmk}\label{Rmk: generalization of the conj}
We do not know of an explicit cyclic group action yielding the cyclic sieving phenomenon 
predicted by Corollary \ref{criterion in type A} in general.
However, in the special case where $n=m$ and $\la=(a^b)$ with $m|a$, 
Rhoades' result \cite[Theorem 1.4]{Rho10} says that 
$ ( \SST_m(\la),  \lan \pr \ran, \allowbreak q^{- \kappa(\la)} s_\la(1,q, \ldots, q^{m-1}) )$ 
exhibits the cyclic sieving phenomenon, where $\pr$ is the promotion operator.
\end{Rmk}

\begin{Ex} Let $\g = \mathfrak{gl}_3$ and $n=4$.  Then we have $\rP^+ = \{ \al_1, \al_2, \al_1+\al_2  \}$ and 
\begin{align*}
\rP^+_1 = \rP^+, \quad \rP^+_2 = \{ \al_1 + \al_2 \}, \quad \rP^+_4 = \emptyset.
\end{align*}
Let $\la = (4)$. It is obvious that $\la$ satisfies the condition $\eqref{Eq: stretched}$. Then we have $\La = 4\epsilon_1$ by $\eqref{eq: la La}$ and 
\begin{align*}
s_\la(x_1,x_2,x_3) =& x_1^4 + x_1^3x_2 + x_1^2x_2^2 + x_1x_2^3 + x_2^4 + x_1^3x_3 + x_1^2x_2x_3 + x_1x_2^2x_3 + x_2^3x_3 \\
& + x_1^2x_3^2 + x_1x_2x_3^2 + x_2^2x_3^2 + x_1x_3^3 + x_2x_3^3 + x_3^4,
\end{align*}
which gives the principal specialization
$$
q^{-\kappa(\la)} s_\la( 1,q, q^2) = q^8 + q^7 + 2q^6 + 2 q^5 + 3 q^4 + 2 q^3 + 2 q^2 + q + 1.
$$
On the other hand,  it follows from $\eqref{Eq: def of a b}$ and Example \ref{Prop: a b} that 
\begin{align*}
\cb_1 &= 1,  \qquad \cb_2 =\frac{4}{2}+1 = 3, \qquad \quad  \cb_4 = \left(\frac{4}{1}+1 \right) \left(\frac{0}{1}+1 \right) \left(\frac{4}{2}+1 \right)=15, \\
\ca_1 &= 1, \qquad \ca_2 = \frac{1}{2} (\cb_2 - \cb_1) = 1, \quad \ca_4 = \frac{1}{4}(\cb_4 - \cb_2) = 3.
\end{align*}
By Theorem \ref{Thm: main}, we have
$$
q^{-\kappa(\la)} s_\la( 1,q, q^2)  = \qdim B(\La) \equiv 1 + (1+q^2) + 3(1+q+q^2+q^3) \pmod{q^4-1}.
$$
Thus, Corollary \ref{criterion in type A} tells us that there exists an action of a cyclic group $C$ of order 4 on $\SST_3(\la)$ such that 
the triple $(\SST_3(\la),  C , q^{-\kappa(\la)} s_\la( 1,q, q^2)  )$ exhibits the cyclic sieving phenomenon.
\end{Ex}

\begin{Ex} Let $\g$ be a simple Lie algebra of type $B_2$.  Then we have 
$$
\cmA = \begin{pmatrix}
2& -1 \\
-2 & 2 
\end{pmatrix}
$$
and $\rP^+ = \{ \al_1,  \al_2,  \al_1+\al_2, \al_1 + 2\al_2  \}$.  Note that $\al_1$ and $\al_1 + 2 \al_2$ are long roots and $\al_2$ and $\al_1 + \al_2$ are short roots.  
We set $n:=2$ and $\Lambda := 2 \varpi_1$.  
Then $(\beta, \Lambda)$ and $\langle \beta^\vee, \Lambda \rangle$ are divisible by $n$ for all $\beta \in \Delta^+$ clearly.

\begin{enumerate}
\item As $(\alpha_1,  \varpi_j) = 2 \delta_{1,j}$ and $ (\alpha_2, \varpi_j) = \delta_{2,j}$,  it follows from \eqref{Eq: q-dim} that 
\begin{align*}
\qdim B(\Lambda) &= \left( \frac{1-q^6}{1-q^2} \right) 
\left( \frac{1-q}{1-q} \right)  \left( \frac{1-q^7}{1-q^3} \right) \left( \frac{1-q^8}{1-q^4} \right) \\
&= 1 + q^2 + q^3 + 2q^4 + q^5 + 2q^6  + q^7 + 2q^8 + q^9 + q^{10} + q^{12}\\
&   \equiv  10+4q   \pmod{q^2-1}.
\end{align*}
Since $ \Delta_1^+ = \Delta^+$ and $\Delta_2^+ = \{  \al_1, \al_1 + 2\al_2 \}$,  it follows from $\eqref{Eq: def of a b}$ and Example \ref{Prop: a b} that 
\begin{align*}
\cb_1 &= \left( \frac{4}{2} + 1 \right) \left( \frac{4} {4} + 1 \right) = 6 ,  \qquad \cb_2 =\left( \frac{4}{2} + 1 \right) \left( \frac{0} {1} + 1 \right) \left( \frac{4}{3} + 1 \right) \left( \frac{4} {4} + 1 \right) =  14 , \\
\ca_1 &= \cb_1 = 6, \qquad\qquad\qquad\quad\qquad \ca_2 = \frac{1}{2} (\cb_2 - \cb_1) = 4.
\end{align*}
By Theorem \ref{Thm: main},  we have
$$
\qdim B(\Lambda)  \equiv  6 + 4(1+q)   \pmod{q^2-1}.
$$

\item As $ \{ \beta^\vee \mid \beta \in \Delta_+ \} = \{ h_1, h_2,  2h_1 +h_2, h_1+h_2  \}$,  it follows from \eqref{Eq: second q-dim} that
\begin{align*}
\qdim^\vee B(\Lambda) &= \left( \frac{1-q^3}{1-q } \right) 
\left( \frac{1-q}{1-q} \right)   \left( \frac{1-q^7}{1-q^3} \right) \left( \frac{1-q^4}{1-q^2} \right) \\
&= 1 + q  + 2q^2 + 2q^3 + 2q^4 + 2q^5  + 2q^6 + q^7 + q^8 \\
&   \equiv  8+6q   \pmod{q^2-1}.
\end{align*}
By Remark \ref{second congruence}, we have 
\begin{align*}
\cb_1^\vee &= \left( \frac{2}{2} + 1 \right)  = 2 ,  \qquad \cb_2^\vee =\left( \frac{2}{1} + 1 \right) \left( \frac{0} {1} + 1 \right) \left( \frac{4}{3} + 1 \right) \left( \frac{2} {2} + 1 \right) =  14 , \\
\ca_1^\vee &= \cb_1^\vee = 2, \qquad\quad\qquad \ca_2^\vee = \frac{1}{2} (\cb_2 - \cb_1) = 6,
\end{align*}
which implies that 
$$
\qdim^\vee B(\Lambda)  \equiv  2 + 6(1+q)   \pmod{q^2-1}.
$$
\end{enumerate}
\end{Ex}

\vskip 2em 

\section{Application to the crystal operator $\cc$ on $\SST_m(\la)$} \label{Application to the crystal operator}

Let $m$ be a positive integer $\ge 2$. Recall that $\SST_m(\la)$ has a $U_q(\mathfrak{gl}_m)$-crystal structure, thus it is equipped with an action of the Weyl group.
Let us consider the operator $\cc := \cs_1\cs_2\cdots \cs_{m-1}$ on $\SST_m(\la)$, where $\cs_i$ is the action on the crystal $\SST_m(\la)$ given by 
the simple reflection $s_i = (i,i+1)$ in the Weyl group. 
Note that the order of $\cc$ is $m$. 
The cyclic action given by this operator was extensively studied in \cite{OP19} in the case where $\ell(\la)<m$ and $\gcd(m, |\la| )=1$.
Under this constraint, it was shown that 
the triple 
$
\left( \SST_m(\la),  \lan \cc \ran, q^{- \kappa(\la)} s_\la(1,q, \ldots, q^{m-1}) \right)
$ 
exhibits the cyclic sieving phenomenon and every orbit is free.
We here focus on the case where $\ell(\la)<m$ and $|\la|$ is divisible by $m$.

Let us collect lemmas which are necessary to develop our arguments.  
Given $T\in \SST_m(\la)$, let ${\rm cont}(T):=(c_1, c_2, \ldots, c_m)$, where $c_i$ is the number of $i$'s occurring in $T$. 

\begin{lem} {\rm (\cite[$\S$2.2. Exercise 2]{F97})}\label{positivity of Kostka number}
Suppose that $\la$ and $\mu$ are partitions of $m$. Then $K_{\la \mu}>0$ if and only if 
$\la \trianglerighteq \mu$, where $K_{\la,\mu}$ is the Kostka number.
\end{lem}

\begin{lem} {\rm (cf. \cite{BS17, HK02})}\label{fixed points}
Let $\la$ be a partition.
The set of fixed points of $\SST_m(\la)$ under the action of $\cc$
is nonempty if and only if $|\la|$ is divisible by $m$, in which case  
it is given by 
\begin{equation}\label{fixed point set}
\left\{T \in \SST_m(\la) : {\rm cont}(T)=\left(\frac {|\la|}{m},\frac {|\la|}{m},\ldots, \frac {|\la|}{m} \right)  \right\}. 
\end{equation}
\end{lem}

\begin{proof}
If $|\la|$ is divisible by $m$, from the definition \eqref{Eq: def of si} of the action of $s_i$ on $\SST_m(\la)$ 
one can infer that  
the set of fixed points of $\SST_m(\la)$ under the action of $\cc$
is given by \eqref{fixed point set}, which is nonempty by Lemma \ref{positivity of Kostka number}. 
Conversely, if $|\la|$ is not divisible by $m$, then
there are no fixed points since ${\rm cont}(T)\ne {\rm cont}(\cc \cdot T)$ for all $T\in \SST_m(\la)$.
\end{proof}

For any two partitions $\la, \mu$, we shall write $\la\sim_m \mu$ if they have the same $m$-core.
For more information on $m$-cores, see \cite[$\S$2.7]{JK81} or \cite[Section I.3. Examples 8]{M95}.
Assume that $\ell(\la) \le m$. 
It is easy to see that if $\la \sim_m 0$, then there exists a unique permutation $w_{\lambda}\in S_m$ such that 
$\la + \delta_m \equiv w_{\lambda}\, \delta_m \pmod m$, where  
$\delta_m=(m-1,m-2,\ldots, 1, 0)$.
The following lemma follows from Examples 17 in \cite[Section I.3]{M95}. 

\begin{lem} \label{congruence of schur poly}
Let $\la$ be a partition of length $\le m$ and let $\zeta_m=e^{2\pi i/m}$.
Then we have   
\begin{equation}\label{modularity of Schur polynomial}
s_\la( 1,q, q^2, \ldots, q^{m-1})|_{q=\zeta_m} = 
\begin{cases}
0 & \text{ if }\la \nsim_m 0,\\
1 & \text{ if $\la \sim_m 0$ and $\epsilon(w_{\lambda})=1$},\\
-1 & \text{ if $\la \sim_m 0$ and $\epsilon(w_{\lambda})=-1$},
\end{cases}	
\end{equation}
where $\epsilon(w)$ is the sign of $w$.
\end{lem}

When $|\la|$ is divisible by $m$,  
we define $\Ula$ to be the semistandard tableau in $\SST_m(\la)$ of content $\mu=\left(\frac {|\la|}{m},\frac {|\la|}{m},\ldots, \frac {|\la|}{m} \right)$
obtained by filling the Young diagram of shape $\la$ with entries in the increasing order 
from left to right and from top to bottom.

\begin{thm}\label{non relatively prime case}
Let $\la$ be a partition of length $<m$.
Assume that there exists a fixed point in $\SST_m(\la)$ under the action of $\cc$, that is, $|\la|$ is divisible by $m$.
Then the following are equivalent.
\begin{enumerate}
\item[(a)] 
The triple $(\SST_m(\la), \lan \cc \ran, \cs_{\la}(1,q,q^2, \ldots, q^{m-1}))$
exhibits the cyclic sieving phenomenon.

\item[(b)]
$\la=(am)$ or $((am)^{m-1})$ for some positive integer $a$.
\end{enumerate}
\end{thm}

\begin{proof}
First, we assume that (a) holds.
In view of \eqref{modularity of Schur polynomial}, one sees that 
there exists only one fixed point, which means that 
$$
\text{ $\Ula$ is the unique fixed point in $ \SST_m(\la)$. } 
$$
Let $\la = (\la_1, \ldots, \la_{\ell})$ with $\la_{\ell} >0$ and let $i_k$ be the first entry of the $k$th row of $\Ula$ for $k=1, \ldots, \ell$. 

Suppose that there is a $k$ such that $| i_{k+1} - i_k | > 1$. Let $b$ be the box $( k+1,1) \in \la$ and $b'$ be the rightmost box in the $k$th row of $\Ula$ whose entry is less than $i_{k+1}$. 
Since $| i_{k+1} - i_k | > 1$, the $k$th row of $\Ula$ contains all $i_k+1$'s. Thus we have $\la_k > |\la|/m$, which says that the entry of the box just below $b'$ is larger than $i_{k+1}$ if it exists.
Setting $T$ to be the tableau obtained from $\Ula$ by swapping the entries of $b$ and $b'$, $T$ is a valid semistandard tableau, which tells us that 
$T$ is also a fixed point. This is a contradiction. Thus we conclude that 
\begin{align} \label{Eq: ik}
i_k = k 
\end{align}
for $k=1, \ldots, \ell$.

If $\ell=1$, then $\la = (a m)$ for some $a$.
Let us first show that $1 < \ell < m-1$ is impossible.
In this case, since $i_\ell=\ell < m-1$, the $\ell$-th row of $\Ula$ should contain all entries equal to $m-1$ and $m$.
Thus $\la_{\ell} > 2  {|\la|} / {m}$ and therefore also $ \la_1 > 2 {|\la|} / m$, contradicting  $ i_2 = 2$. 

Suppose that $\ell=m-1$. 
Since the $\ell$th row of $\Ula$ contains all $m$'s, it follows that  
$\la_{\ell} >  {|\la|} / {m}$ and therefore $\la_k >  {|\la|} / {m}$ for any $k$. 
Combining this inequality with $\eqref{Eq: ik}$, we see that each $k$th row of $\Ula$ has both $k$ and $k+1$.
We assume that there exists an index $k$ such that $\la_{k+1} < \la_k$ and $\la_{j} = \la_k$ for all $ j \le k$.  
Let $b$ be the box $( k, \la_k ) \in \la$ and $b'$ be the leftmost box $(k+1, t)$ of the $(k+1)$st row whose entry is $k+2$.
Note that the entry of $b$ is $k+1$. 
Since $\la_k >  {|\la|} / {m}$, the entry of the box $(k, t)$ is $k$.
Hence the tableaux $T$ obtained from $\Ula$ by swapping the entries of $b$ and $b'$ is a valid semistandard tableau. 
This tells us that 
$T$ is also a fixed point, which is a contradiction. Therefore, $\la$ should be of rectangular shape.
\smallskip 

We now assume that (b) holds.
By \cite[Theorem 1.4]{Rho10}, it suffices to see that our crystal operator $\cc$ coincides with $\pr$.
This is straightforward in the case where $\la=(am)$. So, we assume that $\la=((am)^{m-1})$.  
Pick up any $T\in \SST_m(\la)$.

{\it Case 1.} 
Assume that $1\le i<m$ does not appear in the first column.
Then, for all $i\le j \le m-1$, the $j$th row is filled with only $(j+1)$'s.
Hence, for all $i+1 \le j\le m-1$, both $\sigma_j$ and $\cs_j$ act on $T$ as the identity, where 
$\sigma_j$ is the $j$th \emph{Bender-Knuth involution}.
In case of $1\le j \le i$, ignore all entries not equal to $j$ or $j+1$ and 
all columns that contain both $j$ and $j+1$. 
What remains, which is a sequence of $j$'s immediately followed by $(j+1)$, appears only within one row.
This tells us that both $\sigma_j$ and $\cs_j$ act identically for all $1\le j \le i$. 

{\it Case 2.} 
Assume that $m$ does not appear in the first column. 
In the same manner as above, one sees that both $\sigma_j$ and $\cs_j$ act identically for all $1\le j \le m-1$. 
\end{proof}

\begin{Rmk}
\begin{enumerate}
\item
Let $\la$ be of rectangular shape. Then, as permutations on $\SST_m(\la)$, 
$\cc$ and $\pr$ have the same order, but they are not conjugate in general.
It would be nice to characterize $\la$'s such that $\cc$ and $\pr$ are conjugate, 
equivalently, $\la$'s such that $(\SST_m(\la), \lan \cc \ran, \cs_{\la}(1,q,q^2, \ldots, q^{m-1}))$
exhibits the cyclic sieving phenomenon. 

\item	
Let $\la=(a)$ or $(a^{m-1})$, where $|\lambda|$ is not necessarily divisible by $m$.
Following the proof of the second part in Theorem~\ref{non relatively prime case}, one can also see that 
$\cc$ coincides with $\pr$ as operators on $\SST_m(\la)$.

\item
Let $\g$ be a finite-dimensional simple Lie algebra over $\C$.
The longest Weyl group element $w_0$ defines an involution on the simple roots by $\alpha_i \mapsto \alpha_{i^{\ast}}:=-w_0(\alpha_i).$
Consider the automorphism of $U(\mathfrak g)$ defined by 
$$\phi(e_i)=f_i, \quad \phi(f_i)=e_i, \quad \phi(h_i)=-h_i.$$
Let $ \La^{\vee} := - w_0(\La)$,  $v_\La$ be the highest weight vector and $v_\La^{\rm low}$ the lowest weight vector of $B(\La).$
By \cite[Proposition 21.1.2]{BZ96} and \cite[Proposition 7.1]{L93}, one has the bijection $\phi_{\La}: B(\La) \to B(\La^{\vee})$ satisfying that  $v_\La \mapsto v_{\La^{\vee}}^{\rm low}$ and 
\begin{align*}
\phi_{\La} (\tf_i u)=\te_i \phi_{\La}(u), \quad \phi_{\La}(\te_i u)=\tf_i \phi_{\La}(u), \quad \text{for $u\in B(\La)$ and  $i\in I$.} 
\end{align*}
Using \eqref{Eq: def of si}, it is not difficult to see that  
$\phi_{\La}:B(\La) \to B(\La^{\vee})$ is an isomorphism as $W$-sets.
Hence, in type $A_{m-1}$, we have the isomorphism $\phi_{am\varpi_1}: \SST_m((am)) \buildrel \sim\over \longrightarrow \SST_m((am)^{m-1}) $ as $\lan \cc \ran$-sets. 
This isomorphism $\phi_{am\varpi_1}$ explains why both of $(am)$ and $((am)^{m-1})$ appear in Theorem \ref{non relatively prime case}.
Note that the isomorphism $\phi_{am\varpi_1}$ can be understood as a modification of Sch\"utzenberger's or Luszting's involution.
\end{enumerate}
\end{Rmk}

For each  divisor $d$ of $m$,  
let $ \#\Orb_{\lan \cc \ran}^d(\SST_m(\la))$ denote
the number of orbits of size $d$ in $\SST_m(\la)$ under the action of $\lan \cc \ran$.
If $\la=(am)$ or $((am)^{m-1})$ for any positive integer $a$, then it satisfies the condition \eqref{Assumption}.
This enables us to use \eqref{Eq: def of a b} in computing $ \#\Orb_{\lan \cc \ran}^d(\SST_m(\la))$.

\begin{prop} \label{Prop: formula for orbits}
Assume that $\la$ is either $(am)$ or $((am)^{m-1})$ for any positive integer $a$.
For each  divisor $d$ of $m$, we have 
\begin{equation}\label{formula for number of orbits}
\#\Orb_{\lan \cc \ran}^d(\SST_m(\la))=\frac{1}{d} \sum_{e \mid d} \mu \left(\frac{d}{e} \right)\prod_{1\le k <e}\left(\frac{ae}{k}+1 \right).
\end{equation}
\end{prop}
\begin{proof}
Let $n=m$ and $\La=am \epsilon_1$.
It is not difficult to see that 
\begin{align*}
\rP^+_{d} 
&=\{\alpha \in \rP^+: \text{ ${\rm ht}(\alpha)$ is divisible by $d$} \}\\
&=\bigcup_{1\le k < m/d}\{\alpha\in \rP^+: {\rm ht}(\alpha)=kd\}\\
&=\{\epsilon_i-\epsilon_{i+kd}: 1\le k < m/d \text{ and } 1 \le i \le m-kd\}.
\end{align*}
Since 
$$\frac{ (\epsilon_i-\epsilon_{i+kd}, am\epsilon_1 )}{ (\epsilon_i-\epsilon_{i+kd}, \rho)}=\frac{am}{kd}\delta_{i1},$$
we have 
\begin{equation}\label{explicit computation of b}
b_{d}=\prod_{1\le k <d}\left(\frac{ad}{k}+1 \right).
\end{equation}
For the definition of $b_{d}$, see \eqref{Eq: def of a b}. 
Here, $\delta$ denotes the Kronecker delta function and 
the right hand side of \eqref{explicit computation of b} is understood as $1$ in the case where $d=m$.
Therefore, our assertion follows from \eqref{Eq: def of a b}.
\end{proof}

\begin{Ex}
Note that the right hand side of \eqref{formula for number of orbits} does not depend on the choices of $m$. 
Let $\la=(am)$ or $((am)^{m-1})$.  
For every even positive integer $m$, we have $\#\Orb_{\lan \cc \ran}^2(\SST_m(\la))=a$.
And, for every positive multiple $m$ of $3$, we have $\#\Orb_{\lan \cc \ran}^3(\SST_m(\la))=\frac 32a(a+1)$.
\end{Ex}

In the rest of this section, we assume that $n$ is a prime $p \ge m$ and $\la$ is a partition of length $\le p$. 
Recall that 
\begin{equation}\label{principal specialization}
q^{-\kappa(\la)} \cs_{\la}(1,q,q^2, \ldots, q^{m-1})
=\prod_{1\le i<j \le m} \frac{1-q^{(\la_i -i)-(\la_j -j)}} {1-q^{j-i}}
\end{equation}
(for instance, see \cite[Theorem 7.21.2]{Stan99}).
Let $\mathscr A$ be the set of all partitions $\la$ of length $\le p$ satisfying that 
$\la_i-i \equiv \la_j-j \pmod p$ for some $1 \le i < j \le m$.

\begin{prop} \label{little generalization}
Let $p$ be a prime $\ge m$ and $\la$ a partition of length $\le p$. 
\begin{enumerate}
\item $\cs_{\la}(1,q,q^2, \ldots, q^{m-1}) \equiv 0 \pmod {\Phi_p(q)}$ if and only if $\la \in \mathscr A$.
\item
If $\la \in \mathscr A$, then there exists an action of a cyclic group $C$ of order $p$ on $\SST_m(\la)$ such that 
the triple $(\SST_m(\la), C, q^{-\kappa(\la)} s_\la(1,q, q^2, \ldots, q^{m-1}))$ exhibits the cyclic sieving phenomenon.
\item
There exists an action of a cyclic group $C$ of order $p$ on $\SST_p(\la)$ such that 
the triple $(\SST_p(\la), C, s_\la( 1,q, q^2, \ldots, q^{p-1})$ exhibits the cyclic sieving phenomenon
if and only if 
either $\la \nsim_p 0$ or else $\la \sim_p 0$ and $\epsilon(w_\lambda)=1$.
\end{enumerate}
\end{prop}

\begin{proof}
(1) Note that $\mathbb Q[q]/(\Phi_p(q))$ is a field and $1-q^{j-i}$ appearing in the denominator is a unit for all $1 \le i< j\le n$.
Applying this fact to the right hand side of \eqref{principal specialization}, we obtain the desired result.

(2) Let $\la \in \mathscr A$. By (1), we have that $q^{-\kappa(\la)} s_\la(1,q, q^2, \ldots, q^{m-1})\equiv a_p (1+q+ \cdots + q^{p-1}) \pmod {q^p-1}$ for some positive integer $a_p$. 
Therefore, our assertion can be proven in the same way as in Theorem \ref{Thm: csp and B}~(1).

(3) Due to Lemma \ref{congruence of schur poly}, the condition ``either $\la \nsim_p 0$ or else $\la \sim_p 0$ and $\epsilon(w_\lambda)=1$" is equivalent to saying that 
$s_\la(1,q, q^2, \ldots, q^{m-1})\equiv a_1+ a_p (1+q+ \cdots + q^{p-1}) \pmod {q^p-1}$ for some nonnegative integers $a_1$ and $a_p$. 
Therefore, our assertion can also be proved in the same way as in Theorem \ref{Thm: csp and B}~(1).
\end{proof}

\begin{Rmk}
By virtue of the congruence to Kac \cite[Exercise 10.15]{Kac90}, 
one can derive an analogue of Proposition~\ref{little generalization} for highest weight crystals of any finite type.
\end{Rmk}


\bibliographystyle{amsplain}

\begin{thebibliography}{99}



\bibitem{AA19}
P.~Alexandersson and N.~Amini,
\emph{The cone of cyclic sieving phenomena},
Discrete Math. \textbf{342} (2019), no.~6, 1581--1601.


\bibitem{AOL20}
P.~Alexandersson, E.~K.~O\u{g}uz, and S.~Linusson, \emph{Promotion and cyclic sieving on families of SSYT}, arXiv:2007.10478.



\bibitem{APRU20}
P.~Alexandersson, S.~Pfannerer, M.~Rubey, and J.~Uhlin, \emph{Skew characters and cyclic sieving}, arXiv:2004.01140.



\bibitem{BMS14}
M.~Bennett, B.~Madill, and A.~Stokke,
 \emph{Jeu-de-taquin promotion and a cyclic sieving phenomenon for semistandard hook tableaux},
Discrete Math. \textbf{319} (2014), 62--67.

\bibitem{BST10}
J.~Bandlow, A.~Schilling, and N.~Thi\'ery,
 \emph{On the uniqueness of promotion operators on tensor products of type $A$ crystals},
J. Algebraic Combin. \textbf{31} (2010), no.~2, 217--251.



\bibitem{BZ96}
A.~Berenstein and A.~Zelevinsky, \emph{Canonical bases for the quantum group of type $A_r$ and piecewise-linear combinatorics}, Duke Math. J. \textbf{82} (1996), 
no.~3, 473-–502. 
 

\bibitem{BS17} D.~Bump and A.~Schilling,
\emph{Crystal bases. Representations and combinatorics},
World Scientific Publishing Co. Pte. Ltd., Hackensack, NJ, 2017.


\bibitem{FK14}
B.~Fontaine and J.~Kamnitzer, \emph{Cyclic sieving, rotation, and geometric representation theory}, Selecta Math. (N.S.) \textbf{20} (2014), no.~2, 609--625.


\bibitem{F97}
W.~Fulton, 
\emph{Young tableaux.  
With applications to representation theory and geometry}, London Mathematical Society Student Texts, \textbf{35}. Cambridge University Press, Cambridge, 1997.


\bibitem{HK02} 
J.~Hong and S.-J.~Kang, \emph{Introduction to quantum groups and crystal bases},
Graduate Studies in Mathematics, \textbf{42}. American Mathematical Society, Providence, RI, 2002.

\bibitem{JK81}
G.~James and A.~Kerber, \emph{The representation theory of the symmetric group}, Encyclopedia of Mathematics and its Applications, \textbf{16}. Addison-Wesley Publishing Co., Reading, Mass., 1981.

\bibitem{Kac81}
V.~G.~Kac, \emph{Simple Lie groups and the Legendre symbol}, 
Algebra, Carbondale 1980 (Proc. Conf., Southern Illinois Univ., Carbondale, Ill., 1980), pp. 110–123, 
Lecture Notes in Math., \textbf{848}, Springer, Berlin, 1981. 



\bibitem{Kac90}
V.~G.~Kac, \emph{Infinite-dimensional Lie algebras. Third edition}, Cambridge University Press, Cambridge, 1990. 


\bibitem{Kas90}
M.~Kashiwara, \emph{Crystalizing the q-analogue of universal enveloping algebras},
Comm. Math. Phys. \textbf{133} (1990), no.~2, 249--260.




\bibitem{Kas91}
M.~Kashiwara, \emph{On crystal bases of the $Q$-analogue of universal enveloping algebras},
Duke. Math. J. \textbf{63} (1991), no.~2, 465--516.



\bibitem{Kas93}
M.~Kashiwara, \emph{The crystal base and Littelmann's refined Demazure character formula},
 Duke. Math. J. \textbf{71} (1993), no.~3, 839--858.



\bibitem{LS78}
A.~Lascoux, M.-P.~Sch\"{u}tzenberger, \emph{Le monode plaxique}, Noncommutative structures in algebra and geometric combinatorics (Naples, 1978), pp. 129--156, Quad. ``Ricerca Sci.", \textbf{109}, CNR, Rome, 1981.


\bibitem{L93}
G.~Lusztig, \emph{Introduction to quantum groups}, Progress in Mathematics, \textbf{110}. Birkh\"auser Boston, Inc., Boston, MA, 1993.
 
\bibitem{M95}
I.~G.~Macdonald,
\emph{Symmetric functions and Hall polynomials},
Oxford Mathematical Monographs. Oxford Science Publications. The Clarendon Press, Oxford University Press, New York, 1995.

\bibitem{OP19}
Y.-T.~Oh and E.~Park, \emph{Crystals, semistandard tableaux and cyclic sieving phenomenon}, Electron. J. Combin. \textbf{26} (2019), no.~4, Paper No. 4.39, 19 pp.


\bibitem{RSW04}
V.~Reiner, D.~Stanton, and D.~White, \emph{The cyclic sieving phenomenon}, J. Combin. Theory Ser. A \textbf{108} (2004), no.~1, 17--50.


\bibitem{Rho10}
B.~Rhoades, \emph{Cyclic sieving, promotion, and representation theory}, J. Combin. Theory Ser. A \textbf{117} (2010), no.~1, 38--76.

\bibitem{R2001}
D.~B.~Rush, \emph{Restriction of global bases and Rhoades's theorem}, Adv. Math. \textbf{384} (2021), 107725. 






\bibitem{S11}
B.~Sagan, \emph{The cyclic sieving phenomenon: a survey}, Surveys in combinatorics 2011, 183--233, London Math. Soc. Lecture Note Ser., \textbf{392}, Cambridge Univ. Press, Cambridge, 2011.


\bibitem{Sch72}
M.~P.~Sch\"{u}tzenberger, \emph{Promotion des morphismes d'ensembles ordonn\'es}, Discrete Math. \textbf{2} (1972), 73--94.

\bibitem{Sch77}
M.~P. Sch\"utzenberger, \emph{La correspondance de Robinson}, combinatoire et repr\'esentation du groupe sym\'etrique (Actes Table Ronde CNRS, Univ. Louis-Pasteur Strasbourg, Strasbourg, 1976), pp. 59–113. Lecture Notes in Math., Vol. \textbf{579}, Springer, Berlin, 1977.

\bibitem{Stan99}
R.~Stanley, \emph{Enumerative combinatorics. Vol. 2}, Cambridge Studies in Advanced Mathematics, \textbf{62}. Cambridge University Press, Cambridge, 1999.


\bibitem{Stem94}
J.~R.~Stembridge, \emph{On minuscule representations, plane partitions and involutions in complex Lie groups}, Duke Math. J. \textbf{73} (1994), no.~2, 469--490. 

\bibitem{Stem96}
J.~R.~Stembridge, \emph{Canonical bases and self-evacuating tableaux}, Duke Math. J. \textbf{82} (1996), no.~3, 585--606. 

\bibitem{Wes16}
B.~W.~Westbury, \emph{Invariant tensors and the cyclic sieving phenomenon}, Electron. J. Combin. \textbf{23} (2016), no. 4, Paper No. 4.25, 40 pp.


\end{thebibliography}

\end{document}